# SOME RESULTS ABOUT REVERSES OF CAUCHY-SCHWARZ INEQUALITY IN INNER PRODUCT SPACES*

WANG Gong-bao [1], MA Ji-pu [2]

*(1. Dept. of Basic Courses, Naval University of Engineering, Wuhan 430033, China; 2. Department of Mathematics, Nanjing University, Nanjing 210093, China)*

**Abstract:** Some new reverses of the Cauchy-Schwarz inequality in inner product spaces are presented in this paper. As an application of the main result, a formula for error estimate concerning Cauchy-Schwarz's inequality is provided. The results obtained here complement the recent work of the references.

**MR(2000) Subject Classification:** 46C05, 26D15

## 1. Introduction

Let $H$ be an inner product space over the real or complex number field $K$. The following inequality is known in the literature as Cauchy-Schwarz's inequality:

$$|<x,y>|^2 \le \|x\|^2 \|y\|^2, \quad x, y \in H$$

where $\|u\|^2 = <u,u>, u \in H$. And the equality holds if and only if $x$ and $y$ are linearly dependent.

In recent years, many authors have studied the related topics such as reverses of Cauchy-Schwarz, triangle and Bessel inequalities as well as error estimate problem[1–5]. And a lot of good results regarding the above problems have been obtained. In [6], the author provided a survey on the subject of Cauchy-Schwarz type inequalities. In [7], he gave a new counterpart of Cauchy-Schwarz's inequality and applied it to weighted Hilbert space and weighted sequences space.

In this paper, we'll investigate some new reverses of the Cauchy-Schwarz inequality in real or complex inner product spaces. The results obtained here complement the recent work of the references.

## 2. Main results and proof

Throughout this paper, $H$ denotes an inner product space over the real or complex number field $K$.

Now, the main results of the paper as follows.

**Theorem 1.** Suppose that $x, y \in H, x, y \ne 0$ and $\delta > 0$ satisfy the following condition:

$$x \in B(y, \delta) = \{z \in H \mid \|z - y\| \le \delta\}.$$

(a) If $\|y\| > \delta$, then we have the inequality

$$0 \leq \|x\|^2 \|y\|^2 - |< x, y >|^2 \leq \|x\|^2 \|y\|^2 - [\text{Re} < x, y >]^2 \leq \delta^2 \|x\|^2. \qquad (1)$$

Moreover, the constant $m = 1$ in front of $\delta^2$ is best possible in the sense that it can not be replaced by a smaller positive number.

(b) If $\|y\| = \delta,$ then the following inequality holds.

$$\|x\|^2 \leq 2\,\text{Re} < x, y > \leq 2|< x, y >|. \qquad (2)$$

The constant 2 is best possible in both inequalities.

(c) If $\|y\| < \delta,$ then

$$\|x\|^2 \leq \delta^2 - \|y\|^2 + 2\,\text{Re} < x, y > \leq \delta^2 - \|y\|^2 + 2|< x, y >|, \qquad (3)$$

where the constant 2 is also best possible.

**Proof.** Since $x \in B(y, \delta),$ so $\|x - y\|^2 \leq \delta^2,$ which is equivalent to the following inequality

$$\|x\|^2 + \|y\|^2 - \delta^2 \leq 2\,\text{Re} < x, y >. \qquad (4)$$

(a) If $\|y\| > \delta,$ then we can divide (4) by $\sqrt{\|y\|^2 - \delta^2} > 0$ obtaining

$$\frac{\|x\|^2}{\sqrt{\|y\|^2 - \delta^2}} + \sqrt{\|y\|^2 - \delta^2} \leq \frac{2\,\text{Re} < x, y >}{\sqrt{\|y\|^2 - \delta^2}}. \qquad (5)$$

By using the elementary inequality

$$\lambda a + \frac{1}{\lambda} b \geq 2\sqrt{ab}, \lambda > 0, a, b \geq 0,$$

we have that

$$2\|x\| \leq \frac{\|x\|^2}{\sqrt{\|y\|^2 - \delta^2}} + \sqrt{\|y\|^2 - \delta^2}. \qquad (6)$$

Making use of (5) and (6), we get

$$\|x\|^2 \sqrt{\|y\|^2 - \delta^2} \leq \text{Re} < x, y >. \qquad (7)$$

By taking the square in (7), we obtain the third inequality of (1). The others in (1) are obvious.

To prove the sharpness of the constant $m = 1$, assume that there exists a constant $k > 0$ such that

$$\|x\|^2 \|y\|^2 - [\text{Re} < x, y >]^2 \leq k\delta^2 \|x\|^2, \qquad (8)$$

provided that $x \in B(y, \delta)$ and $\|y\| > \delta.$ We should prove that $k \geq 1$.

In fact, let $\delta = \sqrt{\rho} > 0, \rho \in (0,1), y, z \in H$ with $\|y\| = \|z\| = 1$ and $< y, z >= 0.$ Put $x = y + \sqrt{\rho}z.$ Then $x \in B(y, \delta), \|y\| > \delta$ and $\|x\|^2 = \|y\|^2 + \rho\|z\|^2 = 1 + \rho, \text{Re} < x, y >= \|y\|^2 = 1.$ Therefore, we have

$$\|x\|^2 \|y\|^2 - [\text{Re} <x, y>]^2 = \rho.$$

By using (8), we get that $\rho \leq k\rho(1 + \rho),$ which implies

$$1 \leq k(1 + \rho) \quad \text{for any} \quad \rho \in (0,1). \tag{9}$$

Letting $\rho \to 0^+$, from (9) we get that $k \geq 1$. Hence, the sharpness of the constant in (1) is proved.

(b) Since $\|y\| = \delta,$ it is obvious that the inequality (2) holds. The best constant follows in a similar way to the above proof.

(c) The inequality (3) is obvious from (4). And the best constant can be proved in a similar way to the above. We omit the details.

We have completed the proof of Theorem 1.

Now, by using Theorem 1, we provide the second result in connection to Cauchy-Schwarz's inequality.

**Theorem 2.** Suppose that $x, y \in H, x, y \neq 0$ and $a, b \in K$ satisfy the following condition:

$$\text{Re} <x - ay, by - x> \geq 0. \tag{10}$$

(a) If $\text{Re}(a\bar{b}) > 0,$ then the following inequalities hold.

$$\|x\|^2 \|y\|^2 \leq \frac{1}{4} \frac{(\text{Re}[(\bar{a} + \bar{b}) <x, y>])^2}{\text{Re}(a\bar{b})} \leq \frac{1}{4} \frac{|a+b|^2}{\text{Re}(a\bar{b})} |<x, y>|^2, \tag{11}$$

where the constant $\frac{1}{4}$ is best possible in both inequalities.

(b) If $\text{Re}(a\bar{b}) = 0,$ then we have

$$\|x\|^2 \leq \text{Re}[(\bar{a} + \bar{b}) <x, y>] \leq |a + b| |<x, y>|. \tag{12}$$

(c) If $\text{Re}(a\bar{b}) < 0,$ then

$$\|x\|^2 \leq \text{Re}[(\bar{a} + \bar{b}) <x, y>] - \text{Re}(a\bar{b}) \|y\|^2$$
$$\leq |a+b| |<x, y>| - \text{Re}(a\bar{b}) \|y\|^2. \tag{13}$$

**Proof.** It is not difficult to verify that the condition $\text{Re} <x - ay, by - x> \geq 0$ is equivalent to the following inequality

$$\left\| x - \frac{a+b}{2} y \right\| \leq \frac{1}{2} |a - b| \|y\|. \tag{14}$$

Put $z = \frac{a+b}{2} y$ and $\delta = \frac{1}{2} |a - b| \|y\|$.

Then we have

$$\|z\|^2 - \delta^2 = \frac{|a+b|^2 - |a-b|^2}{4} \|y\|^2 = \text{Re}(a\bar{b}) \|y\|^2. \tag{15}$$

(a) If $\operatorname{Re}(a\overline{b}) > 0,$ then the hypothesis of (a) in Theorem 1 is satisfied . So from the second inequality in (1) we have

$$\|x\|^2\|z\|^2 - [\operatorname{Re} <x,z>]^2 = \frac{1}{4}|a+b|^2\|x\|^2\|y\|^2 - \frac{1}{4}(\operatorname{Re}[(\overline{a}+\overline{b})<x,y>])^2$$

$\leq \frac{1}{4}|a-b|\|x\|^2\|y\|^2.$ It follows that

$$\frac{|a+b|^2 - |a-b|^2}{4}\|x\|^2\|y\|^2 \leq \frac{1}{4}(\operatorname{Re}[(\overline{a}+\overline{b})<x,y>])^2.$$

Hence, from (15) we get $\|x\|^2\|y\|^2 \leq \frac{1}{4}\frac{(\operatorname{Re}[(\overline{a}+\overline{b})<x,y>])^2}{\operatorname{Re}(a\overline{b})}.$ The first inequality in (11) has been proved. The second inequality in (11) is obvious.

To prove the sharpness of the constant $\frac{1}{4}$ in (11), assume that there exists a constant $l > 0$ such that

$$\|x\|^2\|y\|^2 \leq l\frac{(\operatorname{Re}[(\overline{a}+\overline{b})<x,y>])^2}{\operatorname{Re}(a\overline{b})}, \qquad (16)$$

provided that $\operatorname{Re}(a\overline{b}) > 0$ and the inequality in (10) holds. We must prove that $l \geq \frac{1}{4}.$

Put $a,b > 0,$ and let $x = ay.$ Then (10) holds , and from (16) we have

$$a^2\|y\|^4 \leq l\frac{a^2(a+b)^2\|y\|^4}{ab},$$

which implies

$$l \geq \frac{ab}{(a+b)^2} \quad \text{for any} \quad a,b > 0. \qquad (17)$$

It is easy to see that $\sup_{a,b>0}\frac{ab}{(a+b)^2} = \frac{1}{4}.$ It follows from (17) that $l \geq \frac{1}{4}.$ Therefore, the sharpness of the constant $\frac{1}{4}$ in (11) is proved.

(b) and (c) are easy to be proved . We omit the details.
The proof of Theorem 2 is completed.
The following corollary presents a formula for error estimate concerning Cauchy-Schwarz's inequality.

**Corollary.** Under the hypothesis of Theorem 2 and if $\operatorname{Re}(a\overline{b}) > 0,$ then we have the inequality:

$$0 \leq \|x\|^2 \|y\|^2 - |<x,y>|^2 \leq \frac{1}{4}\frac{|a-b|^2}{\operatorname{Re}(a\bar{b})}|<x,y>|^2, \qquad (18)$$

where the constant $\frac{1}{4}$ is best possible in (18).

**Proof.** It is easy to prove the corollary by using (11) in Theorem 2. And the best constant follows in a similar manner to the one in the proof of (a) in Theorem 2. So we omit the details.

It is obvious that the inequality (18) is better than the inequality (1.6) obtained in [8].

*Research supported by the National Natural Science Foundation of China (10271053) and the Science Foundation of Naval University of Engineering (HGDJJ03001).